# A method of finding the asymptotics of q-series based on the convolution of generating functions

Václav Kotěšovec
e-mail: kotesovec2@gmail.com
September 29, 2015

## Abstract

This paper analyzes over 30 types of q-series and the asymptotic behavior of their expansions. A method is described for deriving further asymptotic formulas using convolutions of generating functions with subexponential growth. All variables in the article are integers.

**Theorem 1** (asymptotics of convolution of generating functions with subexponential growth)

Let $r_1 > 0, r_2 > 0, 0 < p < 1$

$$g_1(x) = \sum_{n=0}^{\infty} \alpha_n x^n \qquad \alpha_n \sim v_1 * \frac{exp(r_1\, n^p)}{n^{b_1}}$$

$$g_2(x) = \sum_{n=0}^{\infty} \beta_n x^n \qquad \beta_n \sim v_2 * \frac{exp(r_2\, n^p)}{n^{b_2}}$$

and

$$g(x) = g_1(x) * g_2(x) = \sum_{n=0}^{\infty} a_n x^n$$

Then

$$a_n \sim \frac{\sqrt{2\pi}\, v_1 v_2 \left(r_1^{\frac{1}{1-p}} + r_2^{\frac{1}{1-p}}\right)^{b_1+b_2-\frac{1}{2}(3-p)} exp\left(n^p \left(r_1^{\frac{1}{1-p}} + r_2^{\frac{1}{1-p}}\right)^{1-p}\right)}{\sqrt{(1-p)p}\, r_1^{\frac{b_1-1/2}{1-p}} r_2^{\frac{b_2-1/2}{1-p}} n^{b_1+b_2+\frac{p}{2}-1}}$$

**Proof:**

$$g(x) = g_1(x) * g_2(x) = \sum_{n=0}^{\infty} \alpha_n x^n * \sum_{n=0}^{\infty} \beta_n x^n = \sum_{n=0}^{\infty} \left(\sum_{k=0}^{n} \alpha_{n-k}\, \beta_k\right) x^n$$

We set

$$k = n * m$$

and find the maximum of the term

$$\alpha_{n-k}\,\beta_k \sim v_1 v_2 k^{-b_2}(n-k)^{-b_1} e^{r_2 k^p + r_1(n-k)^p}$$

$$\alpha_{n-k}\,\beta_k \sim v_1 v_2\,(mn)^{-b_2}((1-m)n)^{-b_1}\,e^{r_1\,((1-m)n)^p + r_2\,(mn)^p}$$

The derivative with respect to $m$ is

$$\frac{v_1 v_2\,e^{r_1((1-m)n)^p + r_2(mn)^p}((1-m)(b_2 - pr_2(mn)^p) - b_1 m + mpr_1((1-m)n)^p)}{(m-1)m\,(mn)^{b_2}\,((1-m)n)^{b_1}}$$

Now we solve an equation

$$(1-m)(b_2 - pr_2(mn)^p) - b_1 m + mpr_1((1-m)n)^p = 0$$

But

$$\lim_{n\to\infty} \frac{-b_1 m - b_2 m + b_2 + pr_2 m^{p+1} n^p - pr_2 m^p n^p + mpr_1(1-m)^p n^p}{n^p} = p((m-1)r_2 m^p + mr_1(1-m)^p)$$

and the equation

$$(1-m)\,r_1^{1/(p-1)} = m\,r_2^{1/(p-1)}$$

has a solution

$$m = \frac{1}{1 + r_1^{\frac{1}{1-p}} r_2^{\frac{1}{-1+p}}} = \frac{r_2^{\frac{1}{1-p}}}{r_1^{\frac{1}{1-p}} + r_2^{\frac{1}{1-p}}}$$

We set now

$$k = \left(\frac{r_2^{\frac{1}{1-p}}}{r_1^{\frac{1}{1-p}} + r_2^{\frac{1}{1-p}}}\right) * n + x$$

The value in the maximum is

$$v_1 v_2 \left(\frac{nr_1^{\frac{1}{1-p}}}{r_1^{\frac{1}{1-p}} + r_2^{\frac{1}{1-p}}} - x\right)^{-b_1} \left(\frac{nr_2^{\frac{1}{1-p}}}{r_1^{\frac{1}{1-p}} + r_2^{\frac{1}{1-p}}} + x\right)^{-b_2} \exp\left(r_1\left(\frac{nr_1^{\frac{1}{1-p}}}{r_1^{\frac{1}{1-p}} + r_2^{\frac{1}{1-p}}} - x\right)^p + r_2 \left(\frac{nr_2^{\frac{1}{1-p}}}{r_1^{\frac{1}{1-p}} + r_2^{\frac{1}{1-p}}} + x\right)^p\right)$$

We approximate the exponent with the Maclaurin series (see [2], p. 683 about this method)

$$r_1\left(\frac{nr_1^{\frac{1}{1-p}}}{r_1^{\frac{1}{1-p}} + r_2^{\frac{1}{1-p}}} - x\right)^p + r_2\left(\frac{nr_2^{\frac{1}{1-p}}}{r_1^{\frac{1}{1-p}} + r_2^{\frac{1}{1-p}}} + x\right)^p \sim n^p \left(r_1^{\frac{1}{1-p}} + r_2^{\frac{1}{1-p}}\right)^{1-p}\left(1 - \frac{(1-p)p\,r_1^{\frac{1}{1-p}} r_2^{\frac{1}{p-1}}\left(r_1^{\frac{1}{p-1}} r_2^{\frac{1}{1-p}} + 1\right)^2 x^2}{2n^2}\right)$$

and near $x = 0$ we have

$$a_n \sim \frac{v_1 v_2 \int_{-\infty}^{\infty} \exp\left(n^p \left(r_1^{\frac{1}{1-p}} + r_2^{\frac{1}{1-p}}\right)^{1-p} \left(1 - \frac{(1-p)p\, r_1^{\frac{1}{1-p}} r_2^{\frac{1}{p-1}} \left(r_1^{\frac{1}{p-1}} r_2^{\frac{1}{1-p}} + 1\right)^2 x^2}{2n^2}\right)\right) dx}{\left(\frac{n r_1^{\frac{1}{1-p}}}{r_1^{\frac{1}{1-p}} + r_2^{\frac{1}{1-p}}}\right)^{b_1} \left(\frac{n r_2^{\frac{1}{1-p}}}{r_1^{\frac{1}{1-p}} + r_2^{\frac{1}{1-p}}}\right)^{b_2}}$$

Finally, we obtain

$$a_n \sim \frac{\sqrt{2\pi} v_1 v_2 \, (r_1 r_2)^{\frac{1}{2(1-p)}} \left(r_1^{\frac{1}{1-p}} + r_2^{\frac{1}{1-p}}\right)^{b_1+b_2+\frac{p-3}{2}} \exp\left(n^p \left(r_1^{\frac{1}{1-p}} + r_2^{\frac{1}{1-p}}\right)^{1-p}\right)}{\sqrt{(1-p)p} \, n^{b_1+b_2+\frac{p}{2}-1} \, r_1^{\frac{b_1}{1-p}} r_2^{\frac{b_2}{1-p}}}$$

**QED**.

In program Mathematica we create a functions:

```
convsubexp[v1_, r1_, b1_, v2_, r2_, b2_, p_]:= v1 * v2 *Sqrt[2*Pi] *
(r1^(1/(1-p)) + r2^(1/(1-p)))^(b1 + b2 - (3-p)/2) * E^(n^p * (r1^(1/(1-p)) +
r2^(1/(1-p)))^(1-p)) / (Sqrt[(1-p)*p] * r1^((b1-1/2)/(1-p)) * r2^((b2-1/2)/(1-
p)) * n^(b1 + b2 + p/2 - 1));

convsubexpfun[fun1_, fun2_] := (e1 = PowerExpand[Exponent[fun1, E]]; e2 =
PowerExpand[Exponent[fun2, E]]; en1 = Exponent[fun1, n]; en2 = Exponent[fun2,
n]; ee1 = Exponent[e1, n]; ee2 = Exponent[e2, n];
FullSimplify[convsubexp[fun1/n^en1/Exp[Coefficient[e1, n^ee1]*n^ee1],
Coefficient[e1, n^ee1], -en1, fun2/n^en2/Exp[Coefficient[e2, n^ee2]*n^ee2],
Coefficient[e2, n^ee2], -en2, ee1], n > 0]);
```

Special cases:
$p = 1/2$

$$a_n \sim \frac{2^{3/2} \sqrt{\pi} \, v_1 v_2 \, (r_1^2 + r_2^2)^{b_1+b_2-5/4} \exp\left(\sqrt{(r_1^2 + r_2^2)\, n}\right)}{n^{b_1+b_2-3/4} \, r_1^{2b_1-1} r_2^{2b_2-1}}$$

(This agree with the result by Dewar & Murty, see [20], p. 3 or [21], p. 6)

$p = 2/3$

$$a_n \sim \frac{3\sqrt{\pi} \, v_1 v_2 \, (r_1^3 + r_2^3)^{b_1+b_2-7/6} \exp\left((r_1^3 + r_2^3)^{1/3} n^{2/3}\right)}{r_1^{3b_1-3/2} r_2^{3b_2-3/2} \, n^{b_1+b_2-2/3}}$$

**Corollary**: Self-convolution of generating function with subexponential growth

Let $r > 0, 0 < p < 1$

$$g(x) = \sum_{n=0}^{\infty} \beta_n x^n \qquad \beta_n \sim v * \frac{exp(r * n^p)}{n^b}$$

$$g(x)^2 = \sum_{n=0}^{\infty} a_n x^n$$

we have

$$a_n \sim \frac{v^2 \sqrt{\pi} \; 2^{2b + \frac{p}{2} - 1} \; exp(r \; 2^{1-p} \; n^p)}{\sqrt{pr(1-p)} \; n^{2b + \frac{p}{2} - 1}}$$

Special cases:

$p = 1/2$

$$a_n \sim \frac{v^2 \sqrt{\pi} \; 2^{2b + 1/4} \; exp(r\sqrt{2n})}{\sqrt{r} \; n^{2b - 3/4}}$$

$p = 2/3$

$$a_n \sim \frac{3 \sqrt{\pi} \; v^2 \; 2^{2b - 7/6} \; exp(r \; 2^{1/3} \; n^{2/3})}{\sqrt{r} \; n^{2b - 2/3}}$$

**Theorem 2** (asymptotics of powers of the generating function with subexponential growth)

For $h \geq 1, r > 0, 0 < p < 1$ and

$$g(x) = \sum_{n=0}^{\infty} \beta_n x^n \qquad \beta_n \sim v * \frac{exp(r * n^p)}{n^b}$$

$$g(x)^h = \sum_{n=0}^{\infty} a_n x^n$$

we have an asymptotics

$$a_n \sim v^h \; h^{bh + \frac{hp}{2} - h - \frac{p}{2} + \frac{1}{2}} \left(\frac{2\pi}{(1-p)pr}\right)^{\frac{h-1}{2}} * \frac{exp(r \; h^{1-p} \; n^p)}{n^{bh - \frac{1}{2}(h-1)(2-p)}}$$

**Proof**:

$$v_{h+1} * \frac{exp(r_{h+1} \; n^p)}{n^{b_{h+1}}} = v_h * \frac{exp(r_h \; n^p)}{n^{b_h}} \quad \textbf{convolution} \quad v_1 * \frac{exp(r_1 \; n^p)}{n^{b_1}}$$

We apply Theorem 1

$$v_{h+1} * \frac{exp(r_{h+1} n^p)}{n^{b_{h+1}}} = v_h v_1 * \frac{\sqrt{2\pi} \left(r_h^{\frac{1}{1-p}} + r_1^{\frac{1}{1-p}}\right)^{b_h+b_1+\frac{p-3}{2}} exp\left(n^p \left(r_h^{\frac{1}{1-p}} + r_1^{\frac{1}{1-p}}\right)^{1-p}\right)}{\sqrt{(1-p)p}\, r_1^{\frac{b_1-\frac{1}{2}}{1-p}} r_h^{\frac{b_h-\frac{1}{2}}{1-p}} n^{b_h+b_1+\frac{p}{2}-1}}$$

$$r_{h+1} = \left(r_h^{\frac{1}{1-p}} + r_1^{\frac{1}{1-p}}\right)^{1-p}$$

$$r_h = h^{1-p}\, r$$

$$b_{h+1} = b_h + b_1 + p/2 - 1$$

$$b_h = hb + (p/2 - 1)(h-1)$$

$$v_{h+1} = v_h v_1 * \frac{\sqrt{2\pi} \left(r_h^{\frac{1}{1-p}} + r_1^{\frac{1}{1-p}}\right)^{b_h+b_1+\frac{p-3}{2}}}{\sqrt{(1-p)p}\, r_1^{\frac{b_1-\frac{1}{2}}{1-p}} r_h^{\frac{b_h-\frac{1}{2}}{1-p}}}$$

$$v_{h+1} = v_h * v * \frac{\sqrt{2\pi}\, h^{-bh - \frac{hp}{2} + h + \frac{p}{2} - \frac{1}{2}} (h+1)^{bh + b + \frac{hp}{2} - h - \frac{1}{2}}}{\sqrt{(1-p)pr}}$$

$$v_h = v^h\, h^{bh + \frac{hp}{2} - h - \frac{p}{2} + \frac{1}{2}} \left(\frac{2\pi}{(1-p)pr}\right)^{\frac{h-1}{2}}$$

In program Mathematica:

```
powerself[v_, r_, b_, h_, p_]:= v^h * h^(1/2 - h + b*h - p/2 + h*p/2) *
(2*Pi/((1-p)*p*r))^((h-1)/2) * E^(h^(1-p) * r * n^p) / n^(b*h - (h-1)*(2-
p)/2);

powerselfun[fun_, h_] := (e1 = PowerExpand[Exponent[fun, E]]; en1 =
Exponent[fun, n]; ee1 = Exponent[e1, n]; FullSimplify[powerself[fun/n^en1/
Exp[Coefficient[e1, n^ee1]*n^ee1], Coefficient[e1, n^ee1], -en1, h, ee1], n >
0]);
```

Special cases:

$p = 1/2$

$$a_n \sim \frac{v^h\, 2^{\frac{3(h-1)}{2}}\, \pi^{\frac{h-1}{2}}\, h^{\left(b-\frac{3}{4}\right)h + \frac{1}{4}}\, exp(r\sqrt{hn})}{r^{\frac{h-1}{2}}\, n^{\left(b-\frac{3}{4}\right)h + \frac{3}{4}}}$$

$p = 2/3$

$$a_n \sim \frac{v^h\, 3^{h-1}\, \pi^{\frac{h-1}{2}}\, h^{\left(b-\frac{2}{3}\right)h + \frac{1}{6}}\, exp(r\, h^{1/3}\, n^{2/3})}{r^{\frac{h-1}{2}}\, n^{\left(b-\frac{2}{3}\right)h + \frac{2}{3}}}$$

A very useful is the following theorem, which solve the equation

$$\text{fun0 } \boldsymbol{convolution} \text{ fun1} = \text{fun2}$$

for the given functions fun1 and fun2. All functions fun0, fun1 and fun2 are subexponential.

**Theorem 3**
The solution of the equation

$$v_0 * \frac{exp(r_0\, n^p)}{n^{b_0}} \;\boldsymbol{convolution}\; v_1 * \frac{exp(r_1\, n^p)}{n^{b_1}} = v_2 * \frac{exp(r_2\, n^p)}{n^{b_2}}$$

for $0 < r_1 < r_2$, $0 < p < 1$ is

$$v_0 * \frac{exp(r_0\, n^p)}{n^{b_0}} = \frac{\sqrt{(1-p)p}\, v_2\, r_2^{\frac{1-2b_2}{2(1-p)}} \left(r_2^{\frac{1}{1-p}} - r_1^{\frac{1}{1-p}}\right)^{b_2-b_1-\frac{p}{2}+\frac{1}{2}}}{\sqrt{2\pi}\, v_1\, r_1^{\frac{1-2b_1}{2(1-p)}}} * \frac{exp\left(n^p \left(r_2^{\frac{1}{1-p}} - r_1^{\frac{1}{1-p}}\right)^{1-p}\right)}{n^{b_2-b_1-\frac{p}{2}+1}}$$

**Proof**: We have an equation

$$\frac{\sqrt{2\pi} v_1 v_0\, (r_1 r_0)^{\frac{1}{2(1-p)}} \left(r_1^{\frac{1}{1-p}} + r_0^{\frac{1}{1-p}}\right)^{b_1+b_0+\frac{p-3}{2}} exp\left(n^p \left(r_1^{\frac{1}{1-p}} + r_0^{\frac{1}{1-p}}\right)^{1-p}\right)}{\sqrt{(1-p)p}\, n^{b_1+b_0+\frac{p}{2}-1}\, r_1^{\frac{b_1}{1-p}} r_0^{\frac{b_0}{1-p}}} = v_2 * \frac{exp(r_2\, n^p)}{n^{b_2}}$$

$$\left(r_1^{\frac{1}{1-p}} + r_0^{\frac{1}{1-p}}\right)^{1-p} = r_2$$

$$r_0 = \left(r_2^{\frac{1}{1-p}} - r_1^{\frac{1}{1-p}}\right)^{1-p}$$

$$b_1 + b_0 + p/2 - 1 = b_2$$

$$b_0 = b_2 - b_1 - p/2 + 1$$

$$\frac{\sqrt{2\pi} v_1 v_0\, (r_1 r_0)^{\frac{1}{2(1-p)}} \left(r_1^{\frac{1}{1-p}} + r_0^{\frac{1}{1-p}}\right)^{b_1+b_0+\frac{p-3}{2}}}{\sqrt{(1-p)p}\, r_1^{\frac{b_1}{1-p}} r_0^{\frac{b_0}{1-p}}} = v_2$$

$$v_0 = \frac{\sqrt{(1-p)p}\, v_2\, r_2^{\frac{1-2b_2}{2(1-p)}} \left(r_2^{\frac{1}{1-p}} - r_1^{\frac{1}{1-p}}\right)^{b_2-b_1-\frac{p}{2}+\frac{1}{2}}}{\sqrt{2\pi}\, v_1\, r_1^{\frac{1-2b_1}{2(1-p)}}}$$

In program Mathematica:

```
convsolve0[v1_, r1_, b1_, v2_, r2_, b2_, p_]:= v2 * Sqrt[(1-p)*p] * r2^((1-
2*b2)/(2*(1-p))) * (r2^(1/(1-p)) - r1^(1/(1-p)))^(1/2 - p/2 - b1 + b2) / (v1 *
Sqrt[2*Pi] * r1^((1-2*b1)/(2*(1-p)))) * Exp[(r2^(1/(1-p)) - r1^(1/(1-p)))^(1-
p) * n^p] / n^(1 - p/2 - b1 + b2);

convsolve[fun1_, fun2_] := (e1 = PowerExpand[Exponent[fun1, E]]; e2 =
PowerExpand[Exponent[fun2, E]]; en1 = Exponent[fun1, n]; en2 = Exponent[fun2,
n]; ee1 = Exponent[e1, n]; ee2 = Exponent[e2, n];
FullSimplify[convsolve0[fun1/n^en1/Exp[Coefficient[e1, n^ee1]*n^ee1],
Coefficient[e1, n^ee1], -en1, fun2/n^en2/Exp[Coefficient[e2, n^ee2]*n^ee2],
Coefficient[e2, n^ee2], -en2, ee1], n > 0]);
```

Special cases:

$p = 1/2$

$$v_0 * \frac{exp(r_0\sqrt{n})}{n^{b_0}} = \frac{v_2\, r_2^{1-2b_2}\, (r_2^2 - r_1^2)^{b_2-b_1+1/4}}{2^{3/2}\sqrt{\pi}\, v_1\, r_1^{1-2b_1}} * \frac{exp\left(\sqrt{(r_2^2-r_1^2)n}\right)}{n^{b_2-b_1+3/4}}$$

$p = 2/3$

$$v_0 * \frac{exp(r_0\, n^{2/3})}{n^{b_0}} = \frac{v_2\, r_1^{3b_1-3/2}\, r_2^{3/2-3b_2}\, (r_2^3 - r_1^3)^{b_2-b_1+1/6}}{3\sqrt{\pi}\, v_1} * \frac{exp\left((r_2^3-r_1^3)^{1/3}\, n^{2/3}\right)}{n^{b_2-b_1+2/3}}$$

# Introduction to an asymptotics of q-series

There are several classic results:

Hardy + Ramanujan (1917), see [3], A000041, number of partitions of $n$

$$\prod_{k=1}^{\infty} \frac{1}{1-q^k} \qquad\qquad a_n \sim \frac{exp\left(\pi\sqrt{\frac{2n}{3}}\right)}{4\,n\,\sqrt{3}}$$

---

Meinardus (1954), see [10], p.301, see also [4] Ayoub (1963), A000009, number of partitions of $n$ into distinct parts

$$\prod_{k=1}^{\infty}(1+q^k) \qquad\qquad a_n \sim \frac{exp\left(\pi\sqrt{\frac{n}{3}}\right)}{4*3^{1/4}\,n^{3/4}}$$

---

Ramanujan (1913), A015128, number of overpartitions of $n$

$$\prod_{k=1}^{\infty}\frac{1+q^k}{1-q^k} \qquad\qquad a_n \sim \frac{exp(\pi\sqrt{n})}{8\,n}$$

**convsubexpfun[E^(Pi*Sqrt[2*n/3])/(4*n*Sqrt[3]), E^(Pi*Sqrt[n/3])/(4*3^(1/4)*n^(3/4))]**

---

Following results can be found using Theorem 2 (or with the Meinardus method, see [9])

A000041 (m=1), A000712 (m=2), A000716 (m=3), A023003 (m=4), A144064

$$\prod_{k=1}^{\infty}\frac{1}{(1-q^k)^m} \qquad m>0 \qquad a_n \sim \frac{m^{\frac{m+1}{4}} exp\left(\pi\sqrt{\frac{2mn}{3}}\right)}{2^{\frac{3m+5}{4}} 3^{\frac{m+1}{4}} n^{\frac{m+3}{4}}}$$

**powerselfun[E^(Pi*Sqrt[2*n/3]) / (4*n*Sqrt[3]), m]**

---

A000009 (m=1), A022567 (m=2), A022568 (m=3), A022569 (m=4), A022570 (m=5), ...

$$\prod_{k=1}^{\infty}(1+q^k)^m \qquad m>0 \qquad a_n \sim \frac{m^{1/4} exp\left(\pi\sqrt{\frac{m\,n}{3}}\right)}{2^{\frac{m+3}{2}} 3^{1/4}\, n^{3/4}}$$

**powerselfun[E^(Pi*Sqrt[n/3]) / (4*3^(1/4)*n^(3/4)), m]**

---

A015128 (m=1), A001934 (m=2), A004404 (m=3, alternating)

$$\prod_{k=1}^{\infty}\left(\frac{1+q^k}{1-q^k}\right)^m \qquad m>0 \qquad a_n \sim \frac{m^{\frac{m+1}{4}} exp(\pi\sqrt{m\,n})}{2^{\frac{3(m+1)}{2}}\, n^{\frac{m+3}{4}}}$$

**powerselfun[E^(Sqrt[n]*Pi)/(8*n), m]**

## More general results and Mathematica functions

$$\prod_{k=0}^{\infty} \frac{1}{1-q^{sk+t}} \qquad s>0, t>0, GCD(s,t)=1$$

The following formula found by Ingham (1941), see [5] (cited in [9], p. 394 or in [6], p. 1)

$$a_n \sim \Gamma\left(\frac{t}{s}\right) \pi^{\frac{t}{s}-1} 2^{-\frac{3}{2}-\frac{t}{2s}} 3^{-\frac{t}{2s}} s^{-\frac{1}{2}+\frac{t}{2s}} n^{-\frac{s+t}{2s}} \exp\left(\pi\sqrt{\frac{2n}{3s}}\right)$$

where $\Gamma$ is the Gamma function.

In program Mathematica we create a function

```
partminus[s_, t_]:= Gamma[t/s] * Pi^(t/s-1) * 2^(-3/2-t/(2*s)) * 3^(-t/(2*s))
 * s^(-1/2+t/(2*s)) * n^(-(s+t)/(2*s)) * E^(Pi*Sqrt[2*n/(3*s)]);
```

---

$$\prod_{k=0}^{\infty}(1+q^{sk+t}) \qquad s>0, t>0, GCD(s,t)=1$$

Meinardus (1954), see [10], p.301

$$a_n \sim \frac{\exp\left(\pi\sqrt{\frac{n}{3s}}\right)}{2^{1+\frac{t}{s}} (3s)^{1/4} n^{3/4}}$$

```
partplus[s_, t_]:= E^(Pi*Sqrt[n/(3*s)]) / (2^(1 + t/s) * (3*s)^(1/4) *
n^(3/4));
```

---

Convolution (applied Theorem 1)

$$\prod_{k=0}^{\infty} \frac{1+q^{sk+t}}{1-q^{sk+t}} \qquad s>0, t>0, GCD(s,t)=1$$

$$a_n \sim \frac{\Gamma\left(\frac{t}{s}\right) s^{\frac{t}{2s}-\frac{1}{2}} \pi^{\frac{t}{s}-1} \exp\left(\pi\sqrt{\frac{n}{s}}\right)}{2^{\frac{2t}{s}+1} n^{\frac{t}{2s}+\frac{1}{2}}}$$

```
partratio[s_, t_]:= Gamma[t/s] * s^(t/(2*s) - 1/2) * Pi^(t/s - 1) *
E^(Pi*Sqrt[n/s]) / (2^(2*t/s + 1) * n^(t/(2*s) + 1/2));
```

# Convolutions

$$\prod_{k=0}^{\infty} \frac{1}{(1-q^{sk+t})*(1-q^{ck+d})} \qquad s>0, t>0, c>0, d>0, GCD(s,t,c,d)=1$$

$$a_n \sim \frac{\Gamma\left(\frac{t}{s}\right)\Gamma\left(\frac{d}{c}\right) s^{\frac{t}{2s}-\frac{d}{2c}-\frac{1}{4}} c^{\frac{d}{2c}-\frac{t}{2s}-\frac{1}{4}} (s+c)^{\frac{t}{2s}+\frac{d}{2c}-\frac{1}{4}} \pi^{\frac{t}{s}+\frac{d}{c}-2} \exp\left(\pi\sqrt{\left(\frac{1}{s}+\frac{1}{c}\right)\frac{2n}{3}}\right)}{2^{\frac{t}{2s}+\frac{d}{2c}+\frac{7}{4}} 3^{\frac{t}{2s}+\frac{d}{2c}-\frac{1}{4}} n^{\frac{1}{4}+\frac{t}{2s}+\frac{d}{2c}}}$$

```
convminus[s_, t_, c_, d_]:= Gamma[t/s] * Gamma[d/c] * s^((2*t/s - 2*d/c -
1)/4) * c^((2*d/c - 2*t/s - 1)/4) * (s+c)^((2*t/s + 2*d/c - 1)/4) * Pi^(t/s +
d/c - 2) * E^(Pi*Sqrt[2*(1/s + 1/c)*n/3]) / (2^((2*t/s + 2*d/c + 7)/4) *
3^((2*t/s + 2*d/c - 1)/4) * n^((1 + 2*t/s + 2*d/c)/4));
```

Proof: Applied Theorem 1 for convolution of

$$\prod_{k=0}^{\infty} \frac{1}{(1-q^{sk+t})} * \prod_{k=0}^{\infty} \frac{1}{(1-q^{ck+d})}$$

---

$$\prod_{k=0}^{\infty}(1+q^{sk+t})*(1+q^{ck+d}) \qquad s>0, t>0, c>0, d>0, GCD(s,t,c,d)=1$$

$$a_n \sim \frac{(s+c)^{1/4} \exp\left(\pi\sqrt{\left(\frac{1}{s}+\frac{1}{c}\right)\frac{n}{3}}\right)}{2^{\frac{1}{2}+\frac{t}{s}+\frac{d}{c}} (3sc)^{1/4} n^{3/4}}$$

```
convplus[s_, t_, c_, d_]:= 2^(-1/2 - t/s - d/c) * (s+c)^(1/4) *
E^(Pi*Sqrt[(1/s+1/c)*n/3]) / (3^(1/4) * s^(1/4) * c^(1/4) * n^(3/4));
```

---

$$\prod_{k=0}^{\infty} \frac{1+q^{sk+t}}{1-q^{ck+d}} \qquad s>0, t>0, c>0, d>0, GCD(s,t,c,d)=1$$

$$a_n \sim \frac{c^{\frac{d}{2c}-\frac{1}{2}} (c+2s)^{\frac{d}{2c}} \pi^{\frac{d}{c}-1} \Gamma\left(\frac{d}{c}\right)}{2^{\frac{d}{c}+\frac{s+t}{s}} 3^{\frac{d}{2c}} s^{\frac{d}{2c}} n^{\frac{c+d}{2c}}} \exp\left(\pi\sqrt{\left(\frac{2}{c}+\frac{1}{s}\right)\frac{n}{3}}\right)$$

```
convratio[s_, t_, c_, d_]:= 2^(-d/c - (s+t)/s) * c^(-1/2 + d/(2*c)) * (c +
2*s)^(d/(2*c)) * E^(Sqrt[(2/c + 1/s)*n/3]*Pi) * Pi^(-1 + d/c) * Gamma[d/c] /
(3^(d/(2*c)) * s^(d/(2*c)) * n^((c+d)/(2*c)));
```

## Powers (applied Theorem 2)

$$\prod_{k=0}^{\infty} \frac{1}{(1-q^{sk+t})^m} \qquad m>0, s>0, t>0, GCD(s,t)=1$$

$$a_n \sim \Gamma\left(\frac{t}{s}\right)^m 2^{-\frac{m+5}{4}-\frac{mt}{2s}} 3^{\frac{m-1}{4}-\frac{mt}{2s}} m^{-\frac{m-1}{4}+\frac{mt}{2s}} s^{-\frac{m+1}{4}+\frac{mt}{2s}} \pi^{-m+\frac{mt}{s}} n^{\frac{m-3}{4}-\frac{mt}{2s}} \exp\left(\pi\sqrt{\frac{2mn}{3s}}\right)$$

```
powerminus[s_, t_, m_]:= Gamma[t/s]^m * 2^(-(m+5)/4 - m*t/(2*s)) * 3^((m-1)/4
 - m*t/(2*s)) * m^(-(m-1)/4 + m*t/(2*s)) * s^(-(m+1)/4 + m*t/(2*s)) * Pi^(-m +
 m*t/s) * n^((m-3)/4 - m*t/(2*s)) * E^(Pi*Sqrt[2*m*n/(3*s)]);
```

---

$$\prod_{k=0}^{\infty} (1+q^{sk+t})^m \qquad m>0, s>0, t>0, GCD(s,t)=1$$

$$a_n \sim \frac{2^{\frac{m-3}{2}-\frac{mt}{s}} m^{1/4}}{(3s)^{1/4} n^{3/4}} \exp\left(\pi\sqrt{\frac{mn}{3s}}\right)$$

```
powerplus[s_, t_, m_]:= 2^((m-3)/2 - m*t/s) * m^(1/4) * E^(Pi*Sqrt[m*n/(3*s)])
 / ((3*s)^(1/4) * n^(3/4));
```

---

$$\prod_{k=0}^{\infty} \left(\frac{1+q^{sk+t}}{1-q^{sk+t}}\right)^m \qquad m>0, s>0, t>0, GCD(s,t)=1$$

$$a_n \sim \Gamma\left(\frac{t}{s}\right)^m 2^{\frac{m}{2}-\frac{3}{2}-\frac{2tm}{s}} s^{-\frac{m}{4}-\frac{1}{4}+\frac{tm}{2s}} m^{\frac{1}{4}-\frac{m}{4}+\frac{tm}{2s}} \pi^{\frac{tm}{s}-m} n^{\frac{m}{4}-\frac{3}{4}-\frac{tm}{2s}} \exp\left(\pi\sqrt{\frac{mn}{s}}\right)$$

```
powerratio[s_, t_, m_]:= Gamma[t/s]^m * 2^(m/2 - 3/2 - 2*t*m/s) * s^(-m/4 -
1/4 + t*m/(2*s)) * E^(Pi*Sqrt[m*n/s]) * m^(1/4 - m/4 + t*m/(2*s)) * Pi^(t*m/s
 - m) * n^(m/4 - 3/4 - t*m/(2*s));
```

Special case: $s=2, t=1, m>0$

$$\prod_{k=0}^{\infty} \left(\frac{1+q^{2k+1}}{1-q^{2k+1}}\right)^m$$

A080054 (m=1), A007096 (m=2), A261647 (m=3), A014969 (m=4), A261648 (m=5), A014970 (m=6)

$$a_n \sim \frac{m^{1/4}}{2^{\frac{m}{2}+\frac{7}{4}} n^{3/4}} \exp\left(\pi\sqrt{\frac{mn}{2}}\right)$$

# Various formulas

$$\prod_{k=1}^{\infty} \frac{1+q^k}{1+q^{mk}} \qquad m > 1$$

A000700 (m=2), A003105 (m=3), A070048 (m=4), A096938 (m=5), A261770 (m=6), A097793 (m=7), ...

$$a_n \sim \frac{(m-1)^{1/4}}{2^{3/2}\, 3^{1/4}\, m^{1/4}\, n^{3/4}} \; exp\left(\pi \sqrt{\frac{(m-1)\,n}{3\,m}}\right)$$

```
convplusdenom[m_]:= E^(Pi*Sqrt[(m-1)*n/(3*m)]) * (m-1)^(1/4) / (2^(3/2) *
3^(1/4) * m^(1/4) * n^(3/4));
```

**Proof**:

$$\prod_{k=1}^{\infty} \frac{1+q^k}{1+q^{mk}} = \prod_{k=0}^{\infty} (1+q^{mk+1}) * (1+q^{mk+2}) * \ldots * (1+q^{mk+m-1}) = \prod_{j=1}^{m-1} \prod_{k=0}^{\infty} (1+q^{mk+j})$$

The coefficient of $[q^n]$ in the expansion of product inside is (after the formula by Meinardus with $s = m$ and $t = j$) asymptotic to

$$\frac{exp\left(\pi\sqrt{\frac{n}{3\,m}}\right)}{2^{1+\frac{j}{m}}\, (3\,m)^{1/4}\, n^{3/4}}$$

Now we have a constant (relative to $n$)

$$\prod_{j=1}^{m-1} 2^{1+\frac{j}{m}} = 2^{-\frac{3}{2}(m-1)}$$

and $(m-1)$-fold convolution of

$$\frac{exp\left(\pi\sqrt{\frac{n}{3\,m}}\right)}{(3\,m)^{1/4}\, n^{3/4}}$$

Note that multiple convolution is applied correctly, because $GCD(m, 1, 2, 3, \ldots, m-1) = 1$.
We now set

$$v = 1/(3m)^{1/4},\, r = 1/(3m),\, b = 3/4,\, h = m-1$$

and apply the Theorem 2.

---

$$\prod_{k=1}^{\infty} \left(\frac{1+q^k}{1+q^{mk}}\right)^h \qquad m > 1, h \geq 1$$

Applied Theorem 2

$$a_n \sim \frac{\left(\frac{h(m-1)}{3m}\right)^{1/4} exp\left(\pi\sqrt{\frac{h(m-1)n}{3m}}\right)}{2^{3/2}\, n^{3/4}}$$

$$\prod_{k=1}^{\infty} \frac{1}{(1+q^k)^m} \qquad m > 0$$

A081362 (m=1), A022597 (m=2), A022598 (m=3), A022599 (m=4), ...

$$a_n \sim (-1)^n \frac{m^{1/4}}{2^{7/4} \, 3^{1/4} \, n^{3/4}} \, exp\left(\pi \sqrt{\frac{m\,n}{6}}\right)$$

```
powerplusdenom[m_]:= (-1)^n * E^(Pi*Sqrt[m*n/6]) * m^(1/4) / (2^(7/4) *
3^(1/4) * n^(3/4));
```

**Proof**: From the Euler identity follows

$$\prod_{k=1}^{\infty} \frac{1}{1+q^k} = \prod_{k=1}^{\infty} \frac{1-q^k}{1-q^{2k}} = \prod_{k=1}^{\infty} \frac{(1-q^{2k-1})*(1-q^{2k})}{1-q^{2k}} = \prod_{k=1}^{\infty}(1-q^{2k-1}) = \prod_{k=0}^{\infty}(1-q^{2k+1})$$

$$\prod_{k=1}^{\infty} \left(\frac{1}{1+(-q)^k}\right)^m = \prod_{k=0}^{\infty}(1+q^{2k+1})^m$$

In our notation is then asymptotics of the coefficient of $[q^n]$

$$a_n \sim (-1)^n * \texttt{Simplify[powerplus[2, 1, m]]} = (-1)^n \frac{m^{1/4}}{2^{7/4} \, 3^{1/4} \, n^{3/4}} \, exp\left(\pi \sqrt{\frac{m\,n}{6}}\right)$$

$$\prod_{k=1}^{\infty} \frac{1+q^{mk}}{1+q^k} \qquad m > 1$$

A081360 (m=2), A109389 (m=3), A261734 (m=4), A133563 (m=5), A261736 (m=6), A113297 (m=7), A261735 (m=8), ...

If $m$ is even then

$$a_n \sim (-1)^n \frac{(m+2)^{1/4}}{4\,(6m)^{1/4}\, n^{3/4}} \, exp\left(\pi \sqrt{\frac{(m+2)\,n}{6m}}\right)$$

If $m$ is odd then

$$a_n \sim (-1)^n \frac{(m-1)^{1/4}}{2^{3/2}\,(6m)^{1/4}\, n^{3/4}} \, exp\left(\pi \sqrt{\frac{(m-1)\,n}{6m}}\right)$$

```
convplusnumer[m_]:= If[EvenQ[m], (-1)^n * E^(Pi*Sqrt[(m+2)*n/(6*m)]) *
(m+2)^(1/4) / (4 * (6*m)^(1/4) * n^(3/4)), (-1)^n * E^(Pi*Sqrt[(m-1)*n/(6*m)])
* (m-1)^(1/4) / (2^(3/2) * (6*m)^(1/4) * n^(3/4))];
```

**Proof**: If $m$ is **even** then from the Euler identity follows

$$\prod_{k=1}^{\infty} \frac{1+q^{mk}}{1+q^k} = \prod_{k=1}^{\infty} (1+q^{mk}) * (1-q^{2k-1})$$

$$\prod_{k=1}^{\infty} \frac{1+q^{mk}}{1+(-q)^k} = \prod_{k=0}^{\infty} (1+q^{2k+1}) * \prod_{k=1}^{\infty} (1+q^{mk})$$

and

$$a_n \sim (-1)^n * \texttt{Simplify[convplus[2, 1, m, m]]} = (-1)^n * \frac{\left(\frac{1}{6}+\frac{1}{3m}\right)^{1/4} e^{\pi \sqrt{\frac{1}{6}+\frac{1}{3m}}\sqrt{n}}}{4 n^{3/4}}$$

if $m$ is **odd** then

$$\prod_{k=1}^{\infty} \frac{1+(-q)^{mk}}{1+(-q)^k} = \prod_{k=0}^{\infty} (1+q^{2k+1}) * \prod_{k=0}^{\infty} \left((1-q^{(2k+1)m}) * (1+q^{(2k+2)m})\right)$$

$$\prod_{k=0}^{\infty} \left((1-q^{(2k+1)m}) * (1+q^{(2k+2)m})\right) * \prod_{k=0}^{\infty} \frac{(1+q^{(2k+1)m})}{(1-q^{(2k+1)m})} = \prod_{k=0}^{\infty} (1+q^{(2k+1)m}) * (1+q^{(2k+2)m})$$

```
Simplify[PowerExpand[convsolve[powerratio[2m, m, 1],
convsubexpfun[powerplus[2, 1, 1], convplus[2m, 2m, 2m, m]]]]]
```

$$a_n \sim (-1)^n * \frac{\left(\frac{1}{6}-\frac{1}{6m}\right)^{1/4} e^{\pi \sqrt{\left(\frac{1}{6}-\frac{1}{6m}\right)n}}}{2\sqrt{2}\, n^{3/4}}$$

$$\prod_{k=1}^{\infty} \frac{1-q^{mk}}{1-q^k} \qquad m > 1$$

A000009 (m=2), A000726 (m=3), A001935 (m=4), A035959 (m=5), A219601 (m=6), A035985 (m=7), A261775 (m=8), A104502 (m=9), A261776 (m=10)

The following formula found by Hagis (1971), see [7].[1]

$$a_n \sim \frac{(m-1)^{1/4}}{2^{5/4}\, 3^{1/4}\, m^{3/4}\, n^{3/4}} \exp\left(\pi \sqrt{\frac{2n(m-1)}{3m}}\right)$$

```
hagis[m_]:= E^(Pi*Sqrt[2*n*(m-1)/(3*m)]) * (m-1)^(1/4) / (2 * 6^(1/4) *
m^(3/4) * n^(3/4));
```

---

$$\prod_{k=1}^{\infty} \left(\frac{1-q^{mk}}{1-q^k}\right)^h \qquad m > 1, h \geq 1$$

Applied Theorem 2

$$a_n \sim \frac{h^{1/4}\, (m-1)^{1/4}\, \exp\left(\pi \sqrt{\frac{2h(m-1)n}{3m}}\right)}{2^{5/4}\, 3^{1/4}\, m^{\frac{1}{4}+\frac{h}{2}}\, n^{3/4}}$$

---

$$\prod_{k=1}^{\infty} \frac{1-q^{(2m+1)k}}{1-q^{2k}} \qquad m \geq 1$$

A262346 (m=1), A262364 (m=2)

$$a_n \sim (-1)^n * \frac{(4m+1)^{1/4}}{2^{7/4}\, 3^{1/4}\, (2m+1)^{3/4}\, n^{3/4}} \exp\left(\pi \sqrt{\frac{(4m+1)n}{6(2m+1)}}\right)$$

**Proof**:
We transform the sequence into sequence with nonnegative coefficients using the following identity. If $s$ is **odd** then

$$\prod_{k=1}^{\infty} \left(1-(-q)^{sk}\right) = \prod_{k=1}^{\infty} \left(1+q^{2sk-s}\right) * \left(1-q^{2sk}\right)$$

For $s = 2m+1$ we have

$$\prod_{k=1}^{\infty} \frac{1-(-q)^{(2m+1)k}}{1-q^{2k}} = \prod_{k=1}^{\infty} \frac{\left(1+q^{(2m+1)*(2k-1)}\right)*\left(1-q^{(2m+1)k}\right)*\left(1+q^{(2m+1)k}\right)}{(1-q^k)*(1+q^k)}$$

---

[1] Note that in [8], p.32 is the formula by Hagis cited incorrectly (must be $s \to s-1$ and $24 \to 24\,n$).

$$\prod_{k=1}^{\infty}(1+q^{(2m+1)*(2k-1)}) * \frac{(1-q^{(2m+1)k})}{(1-q^k)} * \frac{(1+q^{(2m+1)k})}{(1+q^k)} * \prod_{k=1}^{\infty}\frac{1+q^k}{1+q^{(2m+1)k}} = \prod_{k=1}^{\infty}(1+q^{(2m+1)*(2k-1)}) * \frac{(1-q^{(2m+1)k})}{(1-q^k)}$$

and the solution (for $m > 0$) follows from

```
Simplify[PowerExpand[convsolve[convplusdenom[2*m+1], convsubexpfun[partplus[4*m+2,
2*m+1], hagis[2*m+1]]]]]
```

Note that for $m = 0$ we have

$$\prod_{k=1}^{\infty}\frac{1-q^{(2m+1)k}}{1-q^{2k}} = \prod_{k=1}^{\infty}\frac{1-q^k}{1-q^{2k}} = \prod_{k=1}^{\infty}\frac{1}{1+q^k}$$

and

$$a_n \sim (-1)^n * \frac{exp\left(\pi\sqrt{\frac{n}{6}}\right)}{2^{7/4}\, 3^{1/4}\, n^{3/4}}$$

---

$$\prod_{k=1}^{\infty}(1-q^k)*(1+q^k)^m \qquad m > 2$$

A085140 (m=3), A261998 (m=4)

$$a_n \sim \frac{exp\left(\pi\sqrt{\frac{(m-2)\,n}{3}}\right)}{2^{\frac{m+1}{2}}\sqrt{n}}$$

**Proof**:

$$\prod_{k=1}^{\infty}(1-q^k)*(1+q^k)^m * \prod_{k=1}^{\infty}\frac{(1+q^k)}{(1-q^k)} = \prod_{k=1}^{\infty}(1+q^k)^{m+1}$$

```
convsolve[powerratio[1, 1, 1], powerplus[1, 1, m+1]]
```

---

$$\prod_{k=1}^{\infty}\frac{1}{(1+q^k)*(1-q^k)^m} \qquad m > 1$$

A002513 (m=2), A029863 (m=3), A262380 (m=4)

$$a_n \sim \frac{(2m-1)^{\frac{m+1}{4}}}{2^{m+1}\, 3^{\frac{m+1}{4}}\, n^{\frac{m+3}{4}}}\, exp\left(\pi\sqrt{\frac{(2m-1)\,n}{3}}\right)$$

**Proof**:
Direct convolution method is not possible, because the sequence with the generating function $\prod_{k=1}^{\infty}\frac{1}{(1+q^k)}$ is alternating. For the correct asymptotics we solve an equation of type fun0 * fun1 = fun2 with the known asymptotics fun1 and fun2 for the **non-alternating** sequences (applied Theorem 3).

$$\prod_{k=1}^{\infty}\frac{1}{(1+q^k)*(1-q^k)^m} * \prod_{k=1}^{\infty}\frac{(1+q^k)}{(1-q^k)} = \prod_{k=1}^{\infty}\frac{1}{(1-q^k)^{m+1}}$$

```
convsolve[powerratio[1, 1, 1], powerminus[1, 1, m+1]]
```

## More examples

The ideal case is if all terms in the numerator are $(1 + q^{c_i k})$ and all terms is the denominator are $(1 - q^{d_j k})$ and $GCD(c_i, d_j) = 1$ for all these coefficients. Lot of sequences can be transformed into such form.

### A100823

$$\prod_{k=1}^{\infty} \frac{(1 + q^k)}{(1 - q^k) * (1 + q^{3k}) * (1 + q^{5k})} = \prod_{k=1}^{\infty} \frac{(1 + q^{5k-1}) * (1 + q^{5k-2}) * (1 + q^{5k-3}) * (1 + q^{5k-4})}{(1 - q^{6k}) * (1 - q^{3k-1}) * (1 - q^{3k-2})}$$

```
FullSimplify[convsubexpfun[convsubexpfun[convplus[5, 1, 5, 4], convplus[5, 2, 5, 3]], convsubexpfun[convminus[3, 1, 3, 2], partminus[6, 6]]]]
```

$$a_n \sim \frac{exp\left(\frac{\pi}{3}\sqrt{\frac{37 n}{5}}\right) \sqrt{37}}{12 \sqrt{5} n}$$

### A147785

$$\prod_{k=1}^{\infty} \frac{(1 - q^{15k})}{(1 - q^{3k}) * (1 - q^{5k})}$$

We solve an equation

$$\prod_{k=1}^{\infty} \frac{(1 - q^{15k})}{(1 - q^{3k}) * (1 - q^{5k})} * \prod_{k=1}^{\infty} \frac{1 + q^k}{1 - q^k} = \prod_{k=1}^{\infty} \frac{(1 - q^{15k})}{(1 - q^k)} * \frac{(1 + q^k)}{(1 - q^{3k}) * (1 - q^{5k})}$$

```
convsolve[partratio[1,1],convsubexpfun[convsubexpfun[hagis[15],partplus[1,1]],
convminus[3,3,5,5]]]
```

$$a_n \sim \sqrt{\frac{7}{5}} * \frac{exp\left(\frac{\pi}{3}\sqrt{\frac{14 n}{5}}\right)}{12 n}$$

# Plane partitions and $k$ in the exponent

A000219 - number of planar partitions of $n$  (MacMahon 1912)

$$\prod_{k=1}^{\infty} \frac{1}{(1-q^k)^k}$$

Wright (1931)[2], see [13]

$$a_n \sim \frac{\zeta(3)^{7/36} \exp\left(3\,\zeta(3)^{1/3} \left(\frac{n}{2}\right)^{2/3} + \frac{1}{12}\right)}{A * 2^{11/36} \sqrt{3\pi}\; n^{25/36}}$$

where $\zeta(3)$ = A002117 is the Riemann Zeta function and A = A074962 is the Glaisher-Kinkelin constant

In the Mathematica
`Zeta[3]^(7/36) * E^(3*Zeta[3]^(1/3) * (n/2)^(2/3) + 1/12) / (Glaisher * Sqrt[3*Pi] * 2^(11/36) * n^(25/36))`

---

A026007 - number of partitions of $n$ into distinct parts, where $n$ different parts of size $n$ are available

$$\prod_{k=1}^{\infty}(1+q^k)^k$$

$$a_n \sim \frac{\zeta(3)^{1/6} \exp\left(\left(\frac{3}{2}\right)^{4/3} \zeta(3)^{1/3}\, n^{2/3}\right)}{2^{3/4}\, 3^{1/3}\, \sqrt{\pi}\, n^{2/3}}$$

`Zeta[3]^(1/6) * E^((3/2)^(4/3) * Zeta[3]^(1/3) * n^(2/3)) / (2^(3/4) * 3^(1/3) * Sqrt[Pi] * n^(2/3))`

---

A156616 (convolution of A000219 and A026007 - applied Theorem 1)

$$\prod_{k=1}^{\infty}\left(\frac{1+q^k}{1-q^k}\right)^k$$

$$a_n \sim \frac{(7\,\zeta(3))^{7/36}}{A * 2^{7/9} \sqrt{3\pi}\, n^{25/36}} \exp\left(\frac{1}{12} + 3 * 2^{-4/3}\, (7\,\zeta(3))^{1/3}\, n^{2/3}\right)$$

`E^(1/12 + 3 * 2^(-4/3) * (7*Zeta[3])^(1/3) * n^(2/3)) * (7*Zeta[3])^(7/36) / (Glaisher * 2^(7/9) * Sqrt[3*Pi] * n^(25/36))`

---

[2] Unfortunately, in many papers is the formula by Wright (see [13]) **cited incorrectly**! For correct version (with $\sqrt{3\pi}$ in the denominator) see [14]. Also in the paper by Almkvist (see [16], p.344), is Wright's formula incomplete, in the denominator should be $\sqrt{3\pi}$, not $\sqrt{\pi}$. In the paper by Steven Finch (see [15]) was this error already corrected.

**Powers** (applied Theorem 2)

$$\prod_{k=1}^{\infty} \frac{1}{(1-q^k)^{mk}} \qquad m > 0$$

A000219 (m=1), A161870 (m=2), A255610 (m=3), A255611 (m=4), A255612 (m=5), A255613 (m=6), A255614 (m=7), A193427 (m=8)

$$a_n \sim \frac{(m\,\zeta(3))^{\frac{m}{36}+\frac{1}{6}} \exp\left(\frac{m}{12} + \frac{3\,(m\,\zeta(3))^{1/3}\,n^{2/3}}{2^{2/3}}\right)}{A^m\,2^{\frac{1}{3}-\frac{m}{36}}\,3^{1/2}\,\pi^{\frac{1}{2}}\,n^{\frac{m}{36}+\frac{2}{3}}}$$

where $\zeta(3)$ = A002117 is the Riemann Zeta function and A = A074962 is the Glaisher-Kinkelin constant

```
powerkminus[m_]:= 2^(m/36 - 1/3) * E^(m/12 + 3 * 2^(-2/3) * m^(1/3) *
Zeta[3]^(1/3) * n^(2/3)) * (m*Zeta[3])^(m/36 + 1/6) / (Glaisher^m * Sqrt[3*Pi]
* n^(m/36 + 2/3));
```

---

$$\prod_{k=1}^{\infty}(1+q^k)^{mk} \qquad m > 0$$

A026007 (m=1), A026011 (m=2), A027346 (m=3), A027906 (m=4)

$$a_n \sim \frac{(m*\zeta(3))^{1/6} \exp\left(\frac{3^{4/3}\,(m\,\zeta(3))^{1/3}\,n^{2/3}}{2^{4/3}}\right)}{2^{\frac{m}{12}+\frac{2}{3}}\,3^{1/3}\,\sqrt{\pi}\,n^{2/3}}$$

```
powerkplus[m_]:= 2^(-2/3 - m/12) * E^((3/2)^(4/3) * m^(1/3) * Zeta[3]^(1/3) *
n^(2/3)) * m^(1/6) * Zeta[3]^(1/6) / (3^(1/3) * Sqrt[Pi] * n^(2/3));
```

---

$$\prod_{k=1}^{\infty}\left(\frac{1+q^k}{1-q^k}\right)^{mk} \qquad m > 0$$

A156616 (m=1), A261386 (m=2), A261389 (m=3)

$$a_n \sim \frac{(7\,m\,\zeta(3))^{\frac{1}{6}+\frac{m}{36}} \exp\left(\frac{m}{12} + \frac{3\,(7\,m\,\zeta(3))^{1/3}\,n^{2/3}}{2^{4/3}}\right)}{A^m\,2^{\frac{2}{3}+\frac{m}{9}}\,\sqrt{3}\,\pi^{\frac{1}{2}}\,n^{\frac{2}{3}+\frac{m}{36}}}$$

```
powerkratio[m_]:= E^(m/12 + 3/2 * (7*m*Zeta[3]/2)^(1/3) * n^(2/3)) * m^(1/6 +
m/36) * (7*Zeta[3])^(1/6 + m/36) / (Glaisher^m * 2^(2/3 + m/9) * Sqrt[3*Pi] *
n^(2/3 + m/36));
```

[A255528](A255528)

$$\prod_{k=1}^{\infty} \frac{1}{(1+q^k)^k}$$

$$a_n \sim (-1)^n * \frac{A * \zeta(3)^{5/36} \, \exp\left(3\,\zeta(3)^{1/3}\, 2^{-5/3}\, n^{2/3} - \frac{1}{12}\right)}{2^{7/9}\, \sqrt{3\pi}\, n^{23/36}}$$

where $\zeta(3)$ = [A002117](A002117) is the Riemann Zeta function and A = [A074962](A074962) is the Glaisher-Kinkelin constant

```
(-1)^n * Glaisher * Zeta[3]^(5/36) * E^(3*Zeta[3]^(1/3)*n^(2/3) / 2^(5/3) - 1/12) / (2^(7/9) * Sqrt[3*Pi] * n^(23/36))
```

**Proof**:

There is an unsigned sequence:

$$\prod_{k=1}^{\infty} \frac{1}{(1+(-q)^k)^k} = \prod_{k=1}^{\infty}(1+q^{2k-1})^{2k-1} * \prod_{k=1}^{\infty}(1-q^{2k})^k = \prod_{k=1}^{\infty}(1+q^{2k-1})^{2k-1} * \prod_{k=1}^{\infty}(1+q^k)^k * \prod_{k=1}^{\infty}(1-q^k)^k$$

$$\prod_{k=1}^{\infty}(1+q^{2k-1})^{2k-1} * \prod_{k=1}^{\infty}(1+q^{2k})^{2k} = \prod_{k=1}^{\infty}(1+q^k)^k$$

Applied Theorem 3 and Theorem 1.

```
ExpandAll[convsolve[powerkminus[1], convsubexpfun[2*convsolve[(powerkplus[2] /. n -> n/2), powerkplus[1]], powerkplus[1]]]]
```

# Meinardus method, case of one simple pole

$$\prod_{k=1}^{\infty} \frac{1}{(1-q^k)^{k^m}} \qquad m > 0$$

A000219 (m=1), A023871 (m=2), A023872 (m=3), A023873 (m=4), A023874 (m=5), A023875 (m=6), A023876 (m=7), A023877 (m=8), A023878 (m=9), A144048

$$a_n \sim \frac{(\Gamma(m+2)\,\zeta(m+2))^{\frac{1-2\zeta(-m)}{2m+4}} \exp\left(\frac{m+2}{m+1}(\Gamma(m+2)\,\zeta(m+2))^{\frac{1}{m+2}} n^{\frac{m+1}{m+2}} + \zeta'(-m)\right)}{\sqrt{2\pi(m+2)}\; n^{\frac{m+3-2\zeta(-m)}{2m+4}}}$$

```
powerkexpminus[m_]:= (Gamma[m+2]*Zeta[m+2])^((1-2*Zeta[-m])/(2*m+4)) * E^((m+2)/(m+1) *
 (Gamma[m+2]*Zeta[m+2])^(1/(m+2)) * n^((m+1)/(m+2)) + Zeta'[-m]) / (Sqrt[2*Pi*(m+2)] *
 n^((m+3-2*Zeta[-m])/(2*m+4)));
```

**Proof**: We use the Meinardus method, for details see [9], [10], [11], [12].

$$\prod_{k=1}^{\infty} \frac{1}{(1-q^k)^{b(k)}}$$

$$b(k) = k^m$$

We have a Dirichlet series

$$d(s) = \sum_{k=1}^{\infty} b(k) * k^{-s} = \sum_{k=1}^{\infty} k^{m-s} = Zeta(s-m)$$

$$d(0) = Zeta(-m)$$
$$d'(0) = Zeta'(-m)$$

In program Mathematica (parameter *r* is a simple pole)

```
meinardusminus[r_]:= (Exp[dd0] * (2*Pi*(r+1))^(-1/2) * (Residue[d[s], {s, r}]
 * Gamma[r+1] * Zeta[r+1])^((1 - 2*d0)/(2*(r+1))) * n^((2*d0 - 2 -
r)/(2*(r+1))) * Exp[n^(r/(r+1))*(1 + 1/r)*(Residue[d[s], {s, r}]*Gamma[r+1] *
Zeta[r+1])^(1/(r+1))]);

b[k_] := k^m;
Clear[d]; d[s_]:= d[s] = Sum[b[k]*k^(-s), {k, 1, Infinity}];
d0 = d[s] /. s -> 0; dd0 = FunctionExpand[(D[d[s], s]) /. s -> 0]; d[s]
```
```
Zeta[-m + s]
```

The function $d(s)$ has a simple pole $s = m + 1$ and we get the result from

```
meinardusminus[m+1]
```

$$\prod_{k=1}^{\infty}(1+q^k)^{k^m} \qquad m>0$$

A026007 (m=1), A027998 (m=2), A248882 (m=3), A248883 (m=4), A248884 (m=5)

$$a_n \sim \frac{2^{\zeta(-m)}\left((1-2^{-m-1})\Gamma(m+2)\zeta(m+2)\right)^{\frac{1}{2m+4}} \exp\left(\frac{m+2}{m+1}\left((1-2^{-m-1})\Gamma(m+2)\zeta(m+2)\right)^{\frac{1}{m+2}} n^{\frac{m+1}{m+2}}\right)}{\sqrt{2\pi(m+2)}\; n^{\frac{m+3}{2m+4}}}$$

```
powerkexpplus[m_]:= 2^(Zeta[-m]) * ((1-2^(-m-1)) * Gamma[m+2] * Zeta[m+2])^(1/(2*m+4))
* E^((m+2)/(m+1) * ((1-2^(-m-1)) * Gamma[m+2] * Zeta[m+2])^(1/(m+2)) * n^((m+1)/(m+2)))
/ (Sqrt[2*Pi*(m+2)] * n^((m+3)/(2*m+4)));
```

**Proof**:
We use the Meinardus method, for details see [9], [10], [11], [12].

$$\prod_{k=1}^{\infty}(1+q^k)^{b(k)}$$

$$b(k)=k^m$$

We have a Dirichlet series

$$d(s)=\sum_{k=1}^{\infty}b(k)*k^{-s}=\sum_{k=1}^{\infty}k^{m-s}=Zeta(s-m)$$

$$d(0)=Zeta(-m)$$

In program Mathematica (parameter *r* is a simple pole)

```
meinardusplus[r_]:= (2^d0*(2*Pi*(r+1))^(-1/2)*(Residue[d[s], {s, r}] *
Gamma[r+1]*(1 - 2^(-r))*Zeta[r+1])^(1/(2*(r+1))) * n^(-(2+r)/(2*(r+1))) *
Exp[n^(r/(r+1))*(1 + 1/r)*(Residue[d[s], {s, r}]*Gamma[r+1]*(1 - 2^(-r)) *
Zeta[r+1])^(1/(r+1))]);

b[k_] := k^m;
Clear[d]; d[s_] := d[s] = Sum[b[k]*k^(-s), {k, 1, Infinity}];
d0 = d[s] /. s -> 0; d[s]
```

Zeta[-m + s]

The function $d(s)$ has a simple pole $s=m+1$ and we get the result from

```
meinardusplus[m+1]
```

Convolution (applied Theorem 1)

$$\prod_{k=1}^{\infty} \left(\frac{1+q^k}{1-q^k}\right)^{k^m} \qquad m \geq 0$$

A156616 (m=1), A206622 (m=2), A206623 (m=3), A206624 (m=4)

$$a_n \sim \frac{\left((2^{m+2}-1)\,\Gamma(m+2)\,\frac{\zeta(m+2)}{2^{2m+3}\,n}\right)^{\frac{1-2\zeta(-m)}{2m+4}} \exp\left(\frac{m+2}{m+1}\left((2^{m+2}-1)\,n^{m+1}\,\Gamma(m+2)\,\frac{\zeta(m+2)}{2^{m+1}}\right)^{\frac{1}{m+2}} + \zeta'(-m)\right)}{\sqrt{(m+2)\,\pi\,n}}$$

```
powerkexpratio[m_]:= ((2^(m+2)-1) * Gamma[m+2] * Zeta[m+2] / (2^(2*m+3) * n))^((1-
2*Zeta[-m])/(2*m+4)) * E^((m+2)/(m+1) * ((2^(m+2)-1) * n^(m+1) * Gamma[m+2] * Zeta[m+2]
/ 2^(m+1))^(1/(m+2)) + Zeta'[-m]) / Sqrt[(m+2)*Pi*n];
```

If $m$ is **even** and $m \geq 2$, then can be simplified as:

$$a_n \sim \frac{\left((2^{m+2}-1)\,\Gamma(m+2)\,\frac{\zeta(m+2)}{2^{2m+3}\,n}\right)^{\frac{1}{2m+4}} \exp\left(\frac{m+2}{m+1}\left((2^{m+2}-1)\,n^{m+1}\,\Gamma(m+2)\,\frac{\zeta(m+2)}{2^{m+1}}\right)^{\frac{1}{m+2}} + (-1)^{m/2}\,\Gamma(m+1)\,\frac{\zeta(m+1)}{2^{m+1}\,\pi^m}\right)}{\sqrt{(m+2)\,\pi\,n}}$$

```
powerkexpratioeven[m_]:= ((2^(m+2)-1) * Gamma[m+2] * Zeta[m+2] / (2^(2*m+3) *
n))^(1/(2*m+4)) * E^((m+2)/(m+1) * ((2^(m+2)-1) * n^(m+1) * Gamma[m+2] * Zeta[m+2] /
2^(m+1))^(1/(m+2)) + (-1)^(m/2) * Gamma[m+1] * Zeta[m+1] / (2^(m+1) * Pi^m)) /
Sqrt[(m+2)*Pi*n];
```

## Meinardus method, case of more poles

$$\prod_{k=1}^{\infty} \frac{1}{(1-q^k)^{mk+c}} \qquad m > 0$$

$$a_n \sim \frac{(m\,\zeta(3))^{\frac{m}{36}+\frac{c}{6}+\frac{1}{6}} \exp\left(\frac{m}{12} - \frac{c^2\,\pi^4}{432\,m\,\zeta(3)} + \frac{c\,\pi^2\,n^{1/3}}{3*2^{4/3}\,(m\,\zeta(3))^{1/3}} + \frac{3\,(m\,\zeta(3))^{1/3}\,n^{2/3}}{2^{2/3}}\right)}{A^m\,2^{\frac{c}{3}+\frac{1}{3}-\frac{m}{36}}\,3^{1/2}\,\pi^{\frac{c+1}{2}}\,n^{\frac{m}{36}+\frac{c}{6}+\frac{2}{3}}}$$

where $\zeta(3)$ = A002117 is the Riemann Zeta function and A = A074962 is the Glaisher-Kinkelin constant

**Proof**:

$$b(k) = mk + c$$

We have a Dirichlet series

$$d(s) = \sum_{k=1}^{\infty} b(k)*k^{-s} = \sum_{k=1}^{\infty} (mk+c)*k^{-s} = m*Zeta(s-1) + c*Zeta(s)$$

I have created a program in the Mathematica.

```
b[k_] := m*k + c;
Clear[d]; d[s_] := d[s] = Sum[b[k]*k^(-s), {k, 1, Infinity}];
d0 = d[s] /. s -> 0; dd0 = FunctionExpand[(D[d[s], s]) /. s -> 0]; d[s]
```

```
m Zeta[-1 + s] + c Zeta[s]
```

Following program has one parameter $r$ = number of poles of $d(s)$, $r \geq 2$. The poles must be a numbers 1, 2, ... , $r$. This is case of equidistant simple poles (see [11], p.21).

```
meinarduspolesminus[r_]:= (If[r == 1, Print["number of poles must be greater than 1"];
Return[0];];
h = r*Residue[d[s], {s, r}] * Gamma[r] * Zeta[r+1]; Clear[ps]; ps[0] = 1;
Do[ps[t] = ps[t] /. Flatten[Solve[
Coefficient[h*Sum[ps[j]*z^j, {j, 0, t}]^(r+1) - Sum[If[i == 0, d0,
i*Residue[d[s], {s, i}] * Gamma[i] * Zeta[i+1]] * h^((r-i)/(r+1)) * z^(r-i)*
Sum[ps[j]*z^j, {j, 0, t}]^(r-i), {i, 0, r}], z^t] == 0, ps[t]]], {t, 1, r+1}];
dn = Expand[h^(1/(r+1))*Sum[ps[j]*z^(j+1), {j, 0, r+1}] /. z->n^(-1/(r+1))];

mm = h^(-d0/(r+1)) * h^((2+r)/(2*(r+1)))/
Sqrt[2*Pi*Residue[d[s], {s, r}] * Gamma[r+2] * Zeta[r+1]]*
n^(-(2 + r - 2*d0)/(2*(r+1))) * Exp[n*dn +
Sum[If[j == 0, dd0, Residue[d[s], {s, j}] * Gamma[j] * Zeta[j+1]]*
Normal[Series[(dn)^(-j), {n, Infinity, 1}]], {j, 0, r}]];

vv = ExpandAll[Simplify[mm, n > 0]]; ee = Exponent[vv, E];
vv/E^ee * E^Sum[If[Exponent[ee[[j]], n]<0, 0, ee[[j]]], {j, 1, Length[ee]}] );
```

The function $d(s)$ has a two poles $s = 1$ and $s = 2$ and we get the result from

```
meinarduspolesminus[2]
```

$$\prod_{k=1}^{\infty}(1+q^k)^{mk+c} \qquad m>0$$

$$a_n \sim \frac{(m*\zeta(3))^{1/6} \exp\left(-\frac{c^2 \pi^4}{1296\, m\, \zeta(3)} + \frac{c\, \pi^2\, n^{1/3}}{2^{5/3}\, 3^{4/3}\, (m\, \zeta(3))^{1/3}} + \frac{3^{4/3}\, (m\, \zeta(3))^{1/3}\, n^{2/3}}{2^{4/3}}\right)}{2^{\frac{m}{12}+\frac{c}{2}+\frac{2}{3}}\, 3^{1/3}\, \sqrt{\pi}\, n^{2/3}}$$

**Proof**:

$$b(k) = mk + c$$

We have a Dirichlet series

$$d(s) = \sum_{k=1}^{\infty} b(k) * k^{-s} = \sum_{k=1}^{\infty}(mk+c) * k^{-s} = m*Zeta(s-1) + c*Zeta(s)$$

```
b[k_] := m*k + c;
Clear[d]; d[s_] := d[s] = Sum[b[k]*k^(-s), {k, 1, Infinity}];
d0 = d[s] /. s -> 0; dd0 = FunctionExpand[(D[d[s], s]) /. s -> 0]; d[s]

m Zeta[-1 + s] + c Zeta[s]
```

Following program has one parameter $r$ = number of poles of $d(s)$, $r \geq 2$.
The poles must be a numbers $1, 2, \ldots, r$.

```
meinarduspolesplus[r_]:= (If[r == 1, Print["number of poles must be greater than 1"];
Return[0];];
h = r*Residue[d[s], {s, r}] * Gamma[r] * (1-2^(-r)) * Zeta[r+1];
Clear[ps]; ps[0] = 1; Do[ps[t] = ps[t] /. Flatten[Solve[
Coefficient[h*Sum[ps[j]*z^j, {j, 0, t}]^(r+1) - Sum[If[i == 0, 0,
i*Residue[d[s], {s, i}] * Gamma[i] * (1-2^(-i))*
Zeta[i+1]]*h^((r-i)/(r+1)) * z^(r-i) *
Sum[ps[j]*z^j, {j, 0, t}]^(r-i), {i, 0, r}], z^t] == 0, ps[t]]], {t, 1, r+1}];
dn = Expand[h^(1/(r+1))*Sum[ps[j]*z^(j+1), {j, 0, r+1}] /. z->n^(-1/(r+1))];

mm = h^((2+r)/(2*(r+1)))/ Sqrt[2*Pi*Residue[d[s], {s, r}]*(1-2^(-r))*
Gamma[r + 2] * Zeta[r+1]] * n^(-(2+r)/(2*(r+1))) * Exp[n*dn +
Sum[If[j == 0, d0*Log[2], Residue[d[s], {s, j}] * Gamma[j]*(1-2^(-j))*
Zeta[j+1]] * Normal[Series[(dn)^(-j), {n, Infinity, 1}]], {j, 0, r}]];

vv = ExpandAll[Simplify[mm, n > 0]]; ee = Exponent[vv, E];
vv/E^ee * E^Sum[If[Exponent[ee[[j]], n]<0, 0, ee[[j]]], {j, 1, Length[ee]}] );
```

The function $d(s)$ has a two poles $s = 1$ and $s = 2$ and we get the result from

```
meinarduspolesplus[2]
```

$$\prod_{k=1}^{\infty}\left(\frac{1+q^k}{1-q^k}\right)^{mk+c} \qquad m > 0$$

$$a_n \sim \frac{(7\,m\,\zeta(3))^{\frac{1}{6}+\frac{c}{6}+\frac{m}{36}} \exp\left(\frac{m}{12} - \frac{c^2\,\pi^4}{336\,m\,\zeta(3)} + \frac{c\,\pi^2\,n^{1/3}}{2^{5/3}\,(7\,m\,\zeta(3))^{1/3}} + \frac{3\,(7\,m\,\zeta(3))^{1/3}\,n^{2/3}}{2^{4/3}}\right)}{A^m\,2^{\frac{2}{3}+\frac{7c}{6}+\frac{m}{9}}\,\sqrt{3}\,\pi^{\frac{c+1}{2}}\,n^{\frac{2}{3}+\frac{c}{6}+\frac{m}{36}}}$$

where $\zeta(3)$ = A002117 is the Riemann Zeta function and A = A074962 is the Glaisher-Kinkelin constant

**Proof:** We apply the following theorem.

**Theorem 4** (modification of Theorem 1 for $p = 2/3$ and an additional term with $exp(n^{1/3})$ )
Let $r_1 > 0, r_2 > 0$

$$g_1(x) = \sum_{n=0}^{\infty} \alpha_n x^n \qquad \alpha_n \sim v_1 * \frac{\exp(s_1\,n^{1/3} + r_1\,n^{2/3})}{n^{b_1}}$$

$$g_2(x) = \sum_{n=0}^{\infty} \beta_n x^n \qquad \beta_n \sim v_2 * \frac{\exp(s_2\,n^{1/3} + r_2\,n^{2/3})}{n^{b_2}}$$

and

$$g(x) = g_1(x) * g_2(x) = \sum_{n=0}^{\infty} a_n x^n$$

Then

$$a_n \sim \frac{3 v_1 v_2 \sqrt{\pi}\,(r_1^3 + r_2^3)^{b_1+b_2-\frac{7}{6}}}{r_1^{3b_1-\frac{3}{2}}\,r_2^{3b_2-\frac{3}{2}}\,n^{b_1+b_2-\frac{2}{3}}} * \exp\left(\frac{(r_2^2 s_1 - r_1^2 s_2)^2}{4 r_1 r_2 (r_1^3 + r_2^3)} + \frac{r_1 s_1 + r_2 s_2}{(r_1^3 + r_2^3)^{1/3}}\,n^{1/3} + (r_1^3 + r_2^3)^{1/3}\,n^{2/3}\right)$$

Proof is same as proof of Theorem 1.

Functions in the Mathematica are

```
convsubexp1323[v1_, s1_, r1_, b1_, v2_, s2_, r2_, b2_]:= 3*v1*v2 * Sqrt[Pi] * ((r1^3 +
r2^3)^(b1 + b2 - 7/6)/(r1^(3*b1 - 3/2) * r2^(3*b2 - 3/2) * n^(b1 + b2 - 2/3))) *
E^((r2^2*s1 - r1^2*s2)^2 / (4*r1*r2*(r1^3 + r2^3)) + (n^(1/3)*(r1*s1 + r2*s2)) / (r1^3
+ r2^3)^(1/3) + n^(2/3)*(r1^3 + r2^3)^(1/3));

convsubexp1323fun[fun1_, fun2_] := (e1 = PowerExpand[Exponent[fun1, E]]; e2 =
PowerExpand[Exponent[fun2, E]]; en1 = Exponent[fun1, n]; en2 = Exponent[fun2, n];
FullSimplify[convsubexp1323[fun1/n^en1/Exp[Coefficient[e1, n^(1/3)]*n^(1/3)]
/Exp[Coefficient[e1, n^(2/3)]*n^(2/3)], Coefficient[e1, n^(1/3)], Coefficient[e1,
n^(2/3)], -en1, fun2/n^en2/Exp[Coefficient[e2, n^(1/3)]*n^(1/3)]/Exp[Coefficient[e2,
n^(2/3)]*n^(2/3)], Coefficient[e2, n^(1/3)], Coefficient[e2, n^(2/3)], -en2], n > 0]);
```

The example is a sequence A261452 ($m = 2, c = -1$)

## Saddle point method

$$\prod_{k=1}^{\infty} \frac{1}{(1-q^k)^{m^k}} \qquad m > 1$$

A034899 (m=2), A144067 (m=3), A144068 (m=4), A144069 (m=5), A144074

$$a_n \sim \frac{m^n \exp\left(2\sqrt{n} - \frac{1}{2} + c_m\right)}{2\sqrt{\pi}\, n^{3/4}}$$

where

$$c_m = \sum_{j=2}^{\infty} \frac{1}{j\,(m^{j-1} - 1)}$$

For a method of proof see [30] or [29].

---

$$\prod_{k=1}^{\infty} (1+q^k)^{m^k} \qquad m > 1$$

A102866 (m=2), A256142 (m=3)

$$a_n \sim \frac{m^n \exp\left(2\sqrt{n} - \frac{1}{2} - c_m\right)}{2\sqrt{\pi}\, n^{3/4}}$$

where

$$c_m = \sum_{j=2}^{\infty} \frac{(-1)^j}{j\,(m^{j-1} - 1)}$$

---

$$\prod_{k=1}^{\infty} \left(\frac{1+q^k}{1-q^k}\right)^{m^k} \qquad m > 1$$

A261519 (m=2), A261520 (m=3)

$$a_n \sim \frac{m^n \exp\left(2\sqrt{2n} - 1 + c\right)}{\sqrt{\pi}\, 2^{3/4}\, n^{3/4}}$$

where

$$c = 2 \sum_{j=1}^{\infty} \frac{1}{(2j+1)(m^{2j} - 1)}$$